\documentclass[10pt]{article}

\usepackage{theorem,amssymb,amsmath,latexsym}

\usepackage{graphicx}

\usepackage{color}                    
\definecolor{red}{rgb}{1,0,.2}        
\definecolor{cjp}{rgb}{.1,.7,.2}        
\definecolor{fmdc}{rgb}{1,0,.8}        

\newcommand{\A}{{\rm A}}
\newcommand{\B}{{\rm B}}
\newcommand{\C}{{\rm C}}
\newcommand{\D}{{\rm D}}
\newcommand{\LP}{L^{\mbox{\tiny\itshape {\rm P}}}}
\newcommand{\BL}{\beta^{\mbox{\tiny\itshape {\rm P}}}}
\newcommand{\GL}{\gamma^{\mbox{\tiny\itshape {\rm P}}}}
\newcommand{\LambdaL}{\Lambda^{\mbox{\tiny\itshape {\rm P}}}}
\newcommand{\GammaL}{\Gamma^{\mbox{\tiny\itshape {\rm P}}}}

\newtheorem{theorem}{Theorem}[section]

\newtheorem{lemma}[theorem]{Lemma}

\newtheorem{proposition}[theorem]{Proposition}
\theorembodyfont{\rm}

\newtheorem{remark}[theorem]{Remark}

\advance \topmargin by -\headheight
\advance \topmargin by -\headsep
\evensidemargin \oddsidemargin
\begin{document}

\begin{centering}
{\Large \textbf{On the representation of the natural numbers by powers of the golden mean}}

\bigskip

{\bf \large Michel Dekking and Ad van Loon}

\bigskip

{\footnotesize  \it Adresses M.~Dekking:  CWI, Amsterdam and Delft University of Technology, Faculty EEMCS.}

\medskip

{\footnotesize \it Email adresses:  Michel.Dekking@cwi.nl, f.m.dekking@tudelft.nl, advloon@upcmail.nl}

\bigskip

{\large \bf Version November 13, 2021}

\end{centering}

\medskip

\begin{abstract}
 \noindent In a base phi representation a natural number is written as a sum of powers of the golden mean $\varphi$.  There are many ways to do this. Well known is the standard representation, introduced by George Bergman in 1957, where a unique representation is obtained by requiring that no consecutive powers $\varphi^n$ and $\varphi^{n+1}$ do occur in the representation. In this paper we introduce a new representation by allowing that the powers $\varphi^0$ and $\varphi^1$ may occur at the same time, but no other consecutive powers. We then argue that this representation is much closer to the classical representation of the natural numbers by powers of an integer than Bergman's standard representation.
\end{abstract}

\medskip


\bigskip


\date{\today}

\bigskip

\section{Introduction}

A natural number $N$ is written in base phi if $N$ has the form
  $$N= \sum_{i=-\infty}^{\infty} a_i \varphi^i,\vspace*{-.0cm}$$
where the $a_i$ are arbitrary non-negative numbers, and where $\varphi:=(1+\sqrt{5})/2$ is the golden mean.

There are many ways to write a number $N$ as a sum of powers of $\varphi$.
In 1957 George Bergman (\cite{Bergman}) proposed restrictions on the digits $a_i$ which entail that the representation becomes unique.
This is generally accepted as \emph{the} representation of the natural numbers in base phi.
A natural number $N$ is written in the Bergman representation  if $N$ has the form
  $$N= \sum_{i=-\infty}^{\infty} d_i \varphi^i,\vspace*{-.0cm}$$
  with digits $d_i=0$ or 1, and where $d_{i+1}d_i = 11$ is not allowed.
  Similarly to base 10 numbers, we  write these representations as
  $$\beta(N) = d_{L}d_{L-1}\dots d_1d_0\cdot d_{-1}d_{-2} \dots d_{R+1}d_R.$$
Here $L$ is the largest positive, and $R$ is the smallest negative power of $\varphi$ that occurs.

\medskip

 The goal of the present paper is to introduce a new representation, which we tendentiously call  the canonical representation, which has properties that are much closer to the classical representation of the natural numbers by powers of an integer than the Bergman representation.
The \emph{canonical representation} of a natural number $N$ by powers of $\varphi$  has the form
  $$N= \sum_{i=-\infty}^{\infty} c_i \varphi^i,\vspace*{-.0cm}$$
  with digits $c_i=0$ or 1, and where $c_{i+1}c_i = 11$ is not allowed, \emph{except} that $c_1c_0 = 11$, as soon as this is possible.
We  write these representations as
  $$\gamma(N) = c_{L}c_{L-1}\dots c_1c_0\cdot c_{-1}c_{-2} \dots c_{R+1}c_R.$$

\medskip


Note that to obtain the canonical representation of $N$, one first looks if there exists a representation of $N$ with $c_1c_0=11$, and no other $c_{i+1}c_i = 11$, and if this is not the case, then $\gamma(N)=\beta(N)$.

\bigskip

The following table compares the two representations. Most of the time, $\gamma(N)=\beta(N)$. The sequence $N=3,7,10,\dots$ for which $\gamma(N)\ne \beta(N)$ is characterized in Proposition \ref{prop: gammabeta}.

\bigskip

\begin{tabular}{|r|c|c|}
\hline
  \; $N^{\phantom{|}}$ & $\beta(N)$ &\; $\gamma(N)$ \\[.0cm]
\hline
1  \; & \;\,1$\cdot$0         &  \:1$\cdot$0               \\
2  \; & \;\,10$\cdot$01       &  \:10$\cdot$01            \\
3  \; & 100$\cdot$01          &  \:11$\cdot$01            \\
4  \; & 101$\cdot$01          &  101$\cdot$01           \\
5  \; &\; 1000$\cdot$1001     & \; 1000$\cdot$1001         \\
6  \; &\; 1010$\cdot$0001     & \; 1010$\cdot$0001         \\
\hline
\end{tabular} \qquad
\begin{tabular}{|r|c|c|}
\hline
  \; $N^{\phantom{|}}$ & $\beta(N)$ &\; $\gamma(N)$ \\[.0cm]
\hline
7  \; & 10000$\cdot$0001    & \; 1011$\cdot$0001         \\
8  \; & 10001$\cdot$0001    &  10001$\cdot$0001        \\
9  \; & 10010$\cdot$0101    &  10010$\cdot$0101        \\
10 \; & 10100$\cdot$0101    &  10011$\cdot$0101        \\
11 \; &10101$\cdot$0101     &  10101$\cdot$0101        \\
12 \; &100000$\cdot$101001  & \, 100000$\cdot$101001     \\
\hline
\end{tabular}

\bigskip

We now come to the heart of the matter. Why does the representation $\gamma(\cdot)$ deserves\footnote{The word `canonical' for our expansion seems to contrast with a free interpretation of Occam's razor: a principle formulated by the 14th-century Franciscan friar William of Ockham. Called Ockam’s razor (often spelled Occam’s razor), it advises you to seek the more economical solution.
 Occam's Razor is the principle that, ``non sunt multiplicanda entia praeter necessitatem" [i.e., ``don't multiply the agents in a theory beyond what's necessary."]”} to be called canonical?
The evidence for this is two-fold. Representations of the natural numbers in number systems can have two important characteristics. These two characteristics might be indicated as `horizontal', and `vertical'. On the one hand, these characteristics are shared by the canonical base phi representation, \emph{and} by the classical base $b$ representation---where, of course, one has to take into account that in base $b$ there are only  digits with non-negative indices. On the other hand, neither one of these characteristics is shared by the Bergman representation.

\section{Addition of base phi representations}

When $N$ and $N'$ are two natural numbers, with base phi representations $a_{L}\dots a_R$ and $a'_{L'}\dots a'_{R'}$, where we allow the digits $a_i, a'_{i}$ to be arbitrary non-negative numbers, then we obtain a base phi representation of $N+N'$, with digits $a_i+a'_i$ for $\max{L,L'}\le i \le \min{R,R'}$, supplementing missing digits by 0's.
In this paper we consider only $\beta(N)+\beta(N')$ and $\gamma(N)+\gamma(N')$. Note that in general $\beta(N)+\beta(N')\ne \beta(N+N')$, and similarly for $\gamma(\cdot)$. Since they represent, nevertheless, the same number, we will write $\beta(N)+\beta(N')\doteq \beta(N+N')$, and  similarly for $\gamma(\cdot)$.

 When we add two numbers in Bergman or canonical representation, then, in general, there is a carry both to the left and (two places) to the right.
 For example
$$\gamma(5) = \gamma(4+1)\doteq \gamma(4)+\gamma(1)= 101\cdot01 + 1\cdot \doteq 102\cdot01 \doteq 110\cdot02 = 1000\cdot1001.$$
Here we used twice that $2\varphi^n=\varphi^{n+1}+\varphi^{n-2}$ for all integers $n$, a direct consequence of $\gamma(2)=10\cdot01.$\,
Note that there is not only a {\it double carry}, but that we also have to get rid of the 11's (except if $c_1c_0 = 11$), by replacing them with 100's.
This is allowed because of the equation $\varphi^{n+2}=  \varphi^{n+1}+\varphi^{n}.$ We call this operation a {\it golden mean shift}.

\section{Existence and uniqueness}

The key to the existence and uniqueness of the canonical representation is the following lemma.

\begin{lemma} \label{lem:gammabeta} A natural number $N$ has a canonical representation $\gamma(N)$ with $c_1c_0=11$ if and only if $N$ has
a Bergman representation $\beta(N)$ with $d_1d_0\cdot d_{-1}=00\cdot0$.
\end{lemma}

\noindent {\it Proof:} The proof is based on the analysis of the Bergman representation  from the paper \cite{Dekk-phi-FQ}.

Let $\beta(N) = d_{L}d_{L-1}\dots d_1d_0\cdot d_{-1}d_{-2} \dots d_{R+1}d_R$. In \cite{Dekk-phi-FQ} the natural numbers $N$ are coded by four letters $\{\A,\B,\C,\D\}$ according to a coding function $T$ as follows:

 We let\\[0.1cm]
 \hspace*{1.1cm} $T(N) =\A$ \;iff\;  $d_1d_0(N)=10$,\quad      $T(N)=\B$  \;iff\; $d_1d_0\cdot d_{-1}(N)=00\cdot 0$,  \\[0.1cm]
 \hspace*{1.1cm} $ T(N)=\C$ \;iff\; $d_0(N)=1$, \hspace*{0.65cm}  $T(N)=\D$   \;iff\;   $d_1d_0\cdot d_{-1}(N)=00\cdot 1$.\\[0.1cm]
This leads to the following scheme.

  \bigskip

\!\!\!\!\begin{tabular}{|r|c|c|}
   \hline
  \; $N^{\phantom{|}}$ & $\beta(N)$ &\!\! $T(N)$\! \\[.0cm]
   \hline
   1\;& \;\!\!\!$1$          & $\C^{\phantom{|}}$ \\
   2\; & \;\:\,\,$10\cdot01$      & $\A$ \\
   3\; & \;$100\cdot01$     & $\B$ \\
   4\; & \;$101\cdot01$     & $\C$ \\
   5\; & \;\:\,$1000\cdot1001$  & $\D$ \\
   6\; & \;\:\,$1010\cdot0001$  & $\A$ \\
   7\; & \;$10000\cdot0001$ & $\B$ \\
   8\; & \;$10001\cdot 0001$  & $\C$\\
   \hline
 \end{tabular}  \quad
 \begin{tabular}{|r|c|c|}
   \hline
   \;$N^{\phantom{|}}$ & $\beta(N)$ & \!\!$T(N)$\!\\[.0cm]
   \hline
   9\; & \;\!\!\!$10010\cdot0101$     &\, $\A^{\phantom{|}}$ \\
   10\; & \;\!\!\!$10100\cdot0101$     & $\B$ \\
   11\; & \;\!\!\!$10101\cdot0101$     & $\C$ \\
   12\; & \;$100000\cdot101001$  & $\D$ \\
   13\; & \;$100010\cdot001001$  & $\A$ \\
   14\; & \;$100100\cdot001001$  & $\B$ \\
   15\; & \;$100101\cdot001001$  & $\C$ \\
   \phantom{a}16\; & \;$101000\cdot100001$  & $\D$ \\
   \hline
 \end{tabular} \quad
 \begin{tabular}{|r|c|c|}
   \hline
   \;$N^{\phantom{|}}$ & $\beta(N)$ & \!\!$T(N)$\! \\[.0cm]
   \hline
   17\; & \;\,\,\,$101010\cdot000001$    & \,$\A^{\phantom{|}}$ \\
   18\; & \;$1000000\cdot000001$     & $\B$ \\
   19\; & \;$1000001\cdot000001$     & $\C$ \\
   20\; & \;$1000010\cdot010001$  & $\A$ \\
   21\; & \;$1000100\cdot010001$  & $\B$ \\
   22\; & \;$1000101\cdot010001$  & $\C$ \\
   23\; & \;$1001000\cdot100101$  & $\D$ \\
   \phantom{a}24\; & \;$1001010\cdot000101$  & $\A$ \\
   \hline
 \end{tabular}

 \vspace*{0.6cm}

 Note that in this table $\A$ is always followed by $\B$, and that $\B$ is always preceded by $\A$. That this is true for all natural numbers $N$ follows directly from Theorem 5.2. in \cite{Dekk-phi-FQ}. With these ingredients we can now give a proof of the lemma.
 The first step is to prove the following claim. Here we use that according to Remark 5.4 in \cite{Dekk-phi-FQ}, Bergman representations with $d_1d_0\cdot d_{-1}(N)=10\cdot 1$ can not occur.

 \medskip

 CLAIM: A natural number $N$ has a canonical representation $\gamma(N)$ with $c_1c_0 = 11$ if and only if $N-1$
has a Bergman representation $\beta(N-1)$ with $d_2d_1d_0\cdot d_{-1} = 010\cdot0$, in other words: $N-1$ has type $\A$.

 \medskip

 [Proof Claim $\Leftarrow$]: Suppose $N-1$ has a Bergman representation $\beta(N-1)$ with $d_2d_1d_0\cdot d_{-1}=010\cdot0$. When we add 1, we find that $\beta(N)\doteq 1\dots 011\cdot0\dots 1$, where $\doteq$ means that we obtain a representation of $N$, but not the Bergman representation. But then clearly we have obtained a representation of $N$ with $c_1c_0=11$, but with no other occurrences of $11$.

 \medskip

 [Proof Claim $\Rightarrow$]: Suppose $\gamma(N)=1\dots 011\cdot 0 \dots 1$.
 When we perform a golden mean shift, we obtain $\gamma(N)\doteq 1\dots 100\cdot 0 \dots 1$.
 Possibly, we have to perform more golden mean shifts in order to obtain a representation of $N$ with no 11. But in any case the result will be of the form  $1\dots 100\cdot 0 \dots 1$.
 By the unicity of the Bergman representation we have found that $\beta(N)= 1\dots 100\cdot 0 \dots 1$. So $N$ is of type $\B$. But then by the fact, given above, $N-1$ must be of type $\A$.

 \medskip

 The lemma now simply follows from the fact, given above,
 that $\A$ is always followed by $\B$ in the $T$-coding of the natural numbers.  \hfill $\Box$

 \medskip

 \begin{proposition} The canonical representation of a natural number is unique.
 \end{proposition}

\noindent {\it Proof:} Suppose $N$ has canonical representations with $c_1c_0\ne 11$. By Lemma \ref{lem:gammabeta}, these representations correspond 1-to-1 to Bergman representations of $N$, so uniqueness follows from the uniqueness of the Bergman representation.

Suppose $N$ has canonical representations with $c_1c_0= 11$. Changing $c_0=1$ to $c_0=0$,  these representations correspond  1-to-1 to Bergman representations of $N\!-\!1$. Again, uniqueness follows from the uniqueness of the Bergman representation. \hfill $\Box$

 \medskip

 How many canonical representation are there in which 11 occurs? It follows from the next proposition that this happens for about 28\% of the natural numbers.

 \begin{proposition} \label{prop: gammabeta} The canonical representation is not equal to the Bergman representation, i.e., $\gamma(N)\ne \beta(N)$, if and only if there exists a natural number $n$, such that $N=\lfloor (\varphi+2)n \rfloor$.
 \end{proposition}

\noindent  {\it Proof:} Since by Lemma \ref{lem:gammabeta}
$\gamma(N)\ne \beta(N)$  if and only if  $N$  is of type $\B,$
 Theorem 5.1 in \cite{Dekk-phi-FQ} gives the result. \hfill $\Box$

\section{The length of representations}\label{sec:length}

In this section we compare the lengths $L+|R|+1$ of the canonical representations $\gamma(N) = c_{L}\dots c_0\cdot c_{-1} \dots c_R$ and the Bergman representations $\beta(N) = d_{L}\dots d_0\cdot d_{-1} \dots d_R.$

Note that in the classical base $b$ representation the natural numbers are partitioned into intervals $B_n:=[b^{n-1},b^{n}-1]$, where the representation of a number $N$ has $n$ digits if and only if $N\in B_n$. For base phi representations, the role of $b^n$ is taken over by the Lucas numbers $L_n$, where $L_0,L_1,L_2,\dots=2,1,3,\dots$ is defined by  $L_0:=2,\, L_1:=1$, and $L_{n}=L_{n-1}+L_{n-2}$ for $n\ge 2$.

It is therefore important to know the representations of the Lucas numbers. The formulas (\ref{eq:gL3}) for $\gamma(L_{2n}+1)$ and (\ref{eq:gL4}) for $\gamma(L_{2n+1}+1)$ will be useful in Section \ref{sec:runs}.

\begin{lemma}\label{lem:Lucas} For all $n\ge 1$ one has
\begin{eqnarray}\label{eq:gLucas}
&\beta(L_{2n})=  10^{2n}\cdot0^{2n-1}1,\quad \gamma(L_{2n}) = [10]^{n-1}11\cdot0^{2n-1}1,\label{eq:gL1}\\
    &\beta(L_{2n+1}) =\gamma(L_{2n+1}) = 1[01]^n\cdot[01]^n,\label{eq:gL2} \\
    &\beta(L_{2n}+1)= \gamma(L_{2n}+1)= 10^{2n-1}1\cdot0^{2n-1}1,\label{eq:gL3}\\
    &\beta(L_{2n+1}+1)=\gamma(L_{2n+1}+1)=  10^{2n+1}\cdot[10]^n01\label{eq:gL4}.
\end{eqnarray}
\end{lemma}

\noindent  {\it Proof:}  The expressions for $\beta(L_{2n})$ and $\beta(L_{2n+1})$ are well-known (see, e.g., \cite{Dekk-phi-FQ}), and easy to prove: they follow directly from $L_{2n}=\varphi^{2n}+\varphi^{-2n}$, and the recursion $L_{2n+1}=L_{2n}+L_{2n-1}$.

When we perform $n$ golden mean shifts on $[10]^{n-1}11\cdot0^{2n-1}1$, we obtain $  10^{2n}\cdot0^{2n-1}1 = \beta(L_{2n})$. This implies the expression for $\gamma(L_{2n})$.

The equality $\gamma(L_{2n+1})=\beta(L_{2n+1})$ follows by an application of Lemma \ref{lem:gammabeta},
since $\beta(L_{2n+1})$ is of type $\C$.

The expression for $\beta(L_{2n}+1)$ and $\gamma(L_{2n}+1)$ follows immediately from Lemma \ref{lem:gammabeta} by adding 1 to the Bergman expansion of $L_{2n}$ in Equation (\ref{eq:gL1}), which yields a valid Bergman expansion for $L_{2n}+1$ which is of type $\C$.

We leave the proof of Equation (\ref{eq:gL4}) to the reader, see also  Lemma 3.3.~(2) in \cite{Hart99}. \hfill $\Box$

\medskip


What are the intervals of constant expansion length for the Bergman representation?\\
As in \cite{Dekk-phi-FQ} we define the so called  \emph{Lucas intervals}\: $\Lambda_{2n}:=[L_{2n},L_{2n+1}]$ and $\Lambda_{2n+1}:=[L_{2n+1}+1, L_{2n+2}-1]$.\\
The next result is Theorem 2.1 in \cite{Hart99}, derived from  Theorem 1 in \cite{Grabner94}.

\begin{proposition}
   The intervals of constant expansion length for the Bergman expansion are the Lucas intervals $\Lambda_n$, $n\ge 1$.
   More precisely: if  $\beta(N) = d_{L}d_{L-1}\dots d_1d_0\cdot d_{-1}d_{-2} \dots d_{R+1}d_R,$
then the left most index $L=L(N)$ and the right most index $R=R(N)$ satisfy\\[-.4cm]
$$L(N)=2n=-R(N) \;{\rm iff}\; N\in \Lambda_{2n}, \quad L(N)=2n\!+1, \;-R(N)=2n+2 \;{\rm iff}\; N\in \Lambda_{2n+1}.$$
 \end{proposition}

\medskip

What is `wrong' with the Bergman Lucas intervals when we compare them with the intervals $B_n$ of constant expansion length for base $b$?
Answer: the odd index intervals are too small compared to the even index intervals: $|\Lambda_{2n}|=L_{2n-1}+1$, and $|\Lambda_{2n+1}|=L_{2n}-1$.

 \medskip

 Our next task is to determine the intervals of constant expansion length for the canonical representation.\\
 We define the \emph{canonical Lucas intervals}
 $$\Gamma_0:=\{1\}, \quad \Gamma_n:=[L_{n}+1,L_{n+1}] \quad {\rm for\;} n\ge 1.$$
 So $\Gamma_1=[2,3],\, \Gamma_2=\{4\},\, \Gamma_3=[5,7],\, \Gamma_4=[8,11]$,\, etc.\\
 Note that $|\Gamma_n|=L_{n+1}-L_n$ for all $n\ge 1$, an expression which is very similar to $|B_n|=b^n-b^{n-1}$ for the classical base $b$ expansion.

 \begin{proposition}\label{prop:length}
   The intervals of constant expansion length for the canonical expansion are the canonical Lucas intervals $\Gamma_n$, $n\ge 1$.
   More precisely: if  $\gamma(N) = c_{L}c_{L-1}\dots c_1c_0\cdot c_{-1}c_{-2} \dots c_{R+1}c_R,$
then the left most index $L=L(N)$ and the right most index $R=R(N)$ satisfy\\[-.4cm]
   $$L(N)=2n=-R(N)  \;{\rm iff}\; N\in \Gamma_{2n}, \quad L(N)=2n\!+1, \;-R(N)=2n+2  \;{\rm iff}\; N\in \Gamma_{2n+1}.$$
 \end{proposition}

 \noindent  {\it Proof:}  Directly from Lemma ~\ref{lem:Lucas} we see that $|\gamma(L_{2n})|= |\beta(L_{2n})|-1$. Therefore we have to move the first number $L_{2n}$ from $\Lambda_{2n}=[L_{2n},L_{2n+1}]$ to  $\Lambda_{2n-1}=[L_{2n-1}+1, L_{2n}-1]$ as a first step to obtain the intervals of constant length expansion for the canonical expansion. This leads exactly to the intervals $\Gamma_n$. It remains to see that this first step is the only change we have to make, i.e., that $|\gamma(N)|= |\beta(N)|$ for all $N\ne L_{2n}$.
 To prove this, note that we can transform the canonical representation to the Bergman representation by a number of golden mean shifts, starting with replacing $011\cdot 0$ in $\gamma(N)$ by $100\cdot 0$. A second golden mean shift will follow if and only if $1011\cdot 0$
 occurs in $\gamma(N)$, and then  $0000\cdot 0$ occurs in $\beta(N)$. Now either this process stops before reaching the left end of $\gamma(N)$ and then $|\gamma(N)|=|\beta(N)|$, or it continues to the left end, and then $\beta(N)=10\dots0\cdot d_{-1}\dots d_R$.
 But, by Lemma \ref{lem:Lucas}, this information suffices  to conclude that $N=L_{2n}$ for some natural number $n$. This follows from the observation in \cite{Dekking-structure-phi}, that in general the $\beta^+\!-$part of an expansion $\beta(N)=\beta^+(N)\cdot\beta^-(N)$ determines $N$, since the $\beta^-\!-$part codes a real number smaller than 1.
  \hfill $\Box$

\section{The recursive structure theorem}\label{sec:ReStTh}

 To obtain recursive relations for the Bergman representation is relatively simple for the intervals $\Lambda_{2n}$, but the intervals $\Lambda_{2n+1}=[L_{2n+1}+1, L_{2n+2}-1]$ have to be divided into three subintervals. These three intervals are\\[-.6cm]
 \begin{align*}
I_n:=&[L_{2n+1}+1,\, L_{2n+1}+L_{2n-2}-1],\\
J_n:=&[L_{2n+1}+L_{2n-2},\, L_{2n+1}+L_{2n-1}],\\
K_n:=&[L_{2n+1}+L_{2n-1}+1,\, L_{2n+2}-1].
\end{align*}


\noindent It will be convenient to use the free group versions of words of 0's and 1's. This means that we will write, for example, $(01)^{-1}0001=1^{-1}001$.
We can then formulate the following result from \cite{Dekk-how-to-add-FQ}.

\begin{theorem}{\bf [Recursive structure theorem for the Bergman representation]}\label{th:recBerg}

\noindent{\,\bf I\;} For all $n\ge 1$ and $k=1,\dots,L_{2n-1}$
one has $ \beta(L_{2n}+k) =  \beta(L_{2n})+ \beta(k) = 10\dots0 \,\beta(k)\, 0\dots 01.$

\noindent{\bf II} For all $n\ge 2$ and $k=1,\dots,L_{2n-2}-1$
\begin{align*}
I_n:&\quad \beta(L_{2n+1}+k) = 1000(10)^{-1}\beta(L_{2n-1}+k)(01)^{-1}1001,\\ K_n:&\quad\beta(L_{2n+1}+L_{2n-1}+k)=1010(10)^{-1}\beta(L_{2n-1}+k)(01)^{-1}0001.
\end{align*}
Moreover, for all $n\ge 2$ and $k=0,\dots,L_{2n-3}$
$$\hspace*{0.7cm}J_n:\quad\beta(L_{2n+1}+L_{2n-2}+k) = 10010(10)^{-1}\beta(L_{2n-2}+k)(01)^{-1}001001.$$
\end{theorem}

\medskip

Since the canonical Lucas intervals are only---literally---marginally different from the Bergman Lucas intervals we can transform Theorem \ref{th:recBerg} into a similar result for the canonical representation.

This time the  interval $\Gamma_{2n+1}=[L_{2n+1}+1, L_{2n+2}]$ has to be divided into the three subintervals \\[-.5cm]
 \begin{align*}
I_n:=&[L_{2n+1}+1,\, L_{2n+1}+L_{2n-2}],\\
J_n:=&[L_{2n+1}+L_{2n-2}+1,\, L_{2n+1}+L_{2n-1}],\\
K_n:=&[L_{2n+1}+L_{2n-1}+1,\, L_{2n+2}].
\end{align*}

The result becomes the following.

\begin{theorem}{\bf [Recursive structure theorem for the canonical representation]}\label{th:recCan}

\noindent{\,\bf I\;} For all $n\ge 1$ and $k=1,\dots,L_{2n-1}$
one has $ \gamma(L_{2n}+k) =   10^{2n}\cdot0^{2n-1}1 + \gamma(k) = 10\dots0 \,\gamma(k)\, 0\dots 01.$

\noindent{\bf II} For all $n\ge 2$ and $k=1,\dots,L_{2n-2}$
\begin{align*}
I_n:&\quad \gamma(L_{2n+1}+k) = 1000(10)^{-1} \gamma(L_{2n-1}+k)(01)^{-1}1001,\\
K_n:&\quad \gamma(L_{2n+1}+L_{2n-1}+k)=1010(10)^{-1} \gamma(L_{2n-1}+k)(01)^{-1}0001.
\end{align*}
Moreover, for all $n\ge 2$ and $k=1,\dots,L_{2n-3}$
$$\hspace*{0.7cm}J_n:\quad \gamma(L_{2n+1}+L_{2n-2}+k) = 10010(10)^{-1} \gamma(L_{2n-2}+k)(01)^{-1}001001.$$
\end{theorem}

\noindent  {\it Proof:}  These statements follow directly from Theorem \ref{th:recBerg}, except that we have to do an extra check for the exceptional numbers $N=L_{2n}$, and for the endpoints of the intervals $I_n$ and $K_n$.

 The numbers $N=L_{2n+2}$ \emph{are} in fact the endpoints of the intervals $K_n$, and the recursion formula above remains valid for $k=L_{2n-2}$ as we can see by an application of Lemma \ref{lem:Lucas}:
$$\gamma(L_{2n+2}) = [10]^{n}11\cdot0^{2n+1}1=10[10]^{n-1}11\cdot0^{2n-1}001=1010(10)^{-1} \gamma(L_{2n})(01)^{-1}0001.$$
The endpoint of $I_n$ is equal to $L_{2n+1}+L_{2n-2}=L_{2n}+L_{2n-1}+L_{2n-2}=2L_{2n}$.
For the recursion formula above to remain valid for $k=L_{2n-2}$ in the $I_n$ case we therefore have to prove that
\begin{align*}
\gamma(2L_{2n})&=1000(10)^{-1} \gamma(L_{2n-1}+L_{2n-2})(01)^{-1}1001\\
   &=1000(10)^{-1} \gamma(L_{2n})(01)^{-1}1001\\
   &=1000(10)^{-1} [10]^{n-1}11\cdot0^{2n-1}1(01)^{-1}1001\\
   &=   1000 [10]^{n-2}11\cdot0^{2n-2}1001.
\end{align*}
We leave the proof of this canonical expansion of $2L_{2n}$ as a (non-trivial) exercise to the reader (Hint: pass to the Bergman expansion, and exploit the unicity of the canonical expansions. See also page 3 of \cite{Dekk-how-to-add-FQ})   \hfill $\Box$

\medskip

\begin{remark}\label{rem:shifts}
  It is important to observe the close relationship between the intervals $I_n, J_n, K_n$ and the canonical Lucas intervals:
  $$I_n=\Gamma_{2n-1}+L_{2n},\quad  J_n=\Gamma_{2n-2}+L_{2n+1},\quad K_n= \Gamma_{2n-1}+L_{2n+1}.$$
  Here we use the notation $A+x=\{a+x:a\in A\}$ for a set of real numbers $A$ and a real number $x$.
\end{remark}

\medskip

We end this section with a typical application of the Recursive Structure Theorem, a lemma that will be useful in the next section.

\begin{lemma}\label{lem:R1001} {\rm a)}\, For all $n\ge 2$ one has $c_{-2n+3}(N)=0$ for all $N$ from $\Gamma_{2n}$.

{\rm b)}\, For all $n\ge 2$ one has $c_{-2n+1}(N)=1$ for the first $L_{2n-1}$ numbers $N$ from $\Gamma_{2n+1}$,  and $c_{-2n+1}(N)=0$ for the last $L_{2n-2}$ numbers $N$ from $\Gamma_{2n+1}$.

\end{lemma}

\noindent  {\it Proof of} a): \;  By the Recursive Structure Theorem, Part I, the expansions in $\Gamma_{2n}=[L_{2n}+1,L_{2n+1}]$ look like the expansions in the interval $[1,L_{2n-1}]$. This interval is a union of $\Gamma_0, \Gamma_1,\dots,\Gamma_{2n-1},\Gamma_{2n-2}$.
Except for the last two, the numbers $N$ from these intervals have canonical expansions which have a right endpoint $R(N)$ with $-R(N)\le 2n-4$.
So automatically we have $c_{-2n+3}(N)=0$ for the $N$ from these intervals. From Proposition \ref{prop:length} we see that both last two intervals $\Gamma_{2n-1}$ and $\Gamma_{2n-2}$ contain only numbers $N$ with $-R(N)=2n-2$. But then the digit $c_{-2n+3}(N)$ directly to the left of $c_{-2n+2}(N)=1$ must be equal to 0.

\medskip

\noindent {\it Proof of} b): \; From Proposition \ref{prop:length} we have $-R(N)=2n+2$ for all   $N\in \Gamma_{2n+1}=I_n \cup J_n\cup K_n.$

Using Remark \ref{rem:shifts} we see that the interval $I_n \cup J_n$  has length $|\Gamma_{2n-1}|+ |\Gamma_{2n-2}|=L_{2n-2}+L_{2n-3}=L_{2n-1}$. We see directly from the Recursive Structure Theorem, Part II, that the expansions of the numbers $N$ in $I_n \cup J_n$ have $c_{R+3}(N)c_{R+2}(N)c_{R+1}(N)c_{R}(N)= 1001$.
Here  $R(N)+3=-2n-2+3=-2n+1$. On the other hand we see that the expansions in the interval $K_n$ all have in $c_{R+3}(N)c_{R+2}(N)c_{R+1}(N)c_{R}(N)=0001$. These two observations  imply part b).    \hfill $\Box$

\section{Vertical runs}\label{sec:runs}

When we make a table of the classical base $b$ expansions of the natural numbers one observes a very regular structure of the runs of the digits in the columns of the table. As an example, consider the case $b=2$ of the binary expansion.

  \bigskip

\!\!\!\!\begin{tabular}{|r|r|}
\hline
\; $N^{\phantom{|}}$ & expansion  \\[.0cm]
\hline
   0\; &   $0$    \\
   1\; &   $1$     \\
   2\; &  $10$    \\
   3\; &  $11$    \\
   4\; &  $100$   \\
   5\; &  $101$   \\
   6\; &  $110$   \\
   7\; &  $111$   \\
\hline
\end{tabular} \quad
\begin{tabular}{|r|r|}
\hline
\; $N^{\phantom{|}}$ & expansion  \\[.0cm]
   \hline
   8\; &  $1000$  \\
   9\;&   $1001$     \\
   10\; &  $1010$    \\
   11\; &  $1011$    \\
   12\; &  $1100$   \\
   13\; &  $1101$   \\
   14\; &  $1110$   \\
   15\; &  $1111$   \\
   \hline
\end{tabular} \quad
\begin{tabular}{|r|r|}
\hline
\; $N^{\phantom{|}}$ & expansion  \\[.0cm]
   \hline
   16\; &  $10000$  \\
   17\;&   $10001$    \\
   18\; &  $10010$    \\
   19\; &  $10011$    \\
   20\; &  $10100$   \\
   21\; &  $10101$   \\
   22\; &  $10110$   \\
   23\; &  $10111$   \\
   \hline
\end{tabular}

\medskip

In digit position $i$, for $i\ge 0$,  only runs of $2^i$ 1's occur---separated by runs of $2^i$ 0's.

\medskip

For the Bergman expansion there is no such  regularity: vertical runs of 1's of length 1,2,3,4,5,6 and 7 do occur.
This is completely different for the canonical expansion: see Theorem \ref{th:runs}. Part of the proof of this theorem is provided by the following lemma, that compares digits of the numbers at the end and the beginning of canonical Lucas intervals.

\begin{lemma} \label{lem:ortho} Let $N$ have canonical expansion $\gamma(N)=c_L(N)\dots c_R(N)$,
where we add 0's when comparing two expansions, for example $c_{2n+1}(L_{2n})=0$.
   Then for all $n\ge1$:
\begin{align*}
 {\rm [OE]}\qquad &   c_i(L_{2n})=1 \;\,and\;\, c_i(L_{2n}+1)=1 \qquad {\rm happens\:if\:and\:only\:if} \;i=0, \;{\rm or}\; i= -2n, \\
 {\rm [EO]}\qquad &    c_i(L_{2n+1})=1 \;and\; c_i(L_{2n+1}+1)=1 \quad \;{\rm does\: not\: happen\: for\: any} \;i=-2n-2,\dots, 2n+2.
\end{align*}
\end{lemma}

\noindent  {\it Proof:} This follows directly from Lemma \ref{lem:Lucas}. See, e.g., Equation (\ref{eq:gL1}) and Equation (\ref{eq:gL3}) for [OE].    \hfill $\Box$

\medskip

In the Lemma, [OE] refers to the indices of the two successive Lucas intervals  $\Gamma_{2n-1}$, $\Gamma_{2n}$,  [EO] to the indices of the two successive Lucas intervals $\Gamma_{2n}$, $\Gamma_{2n+1}$.

\begin{theorem}\label{th:runs}  In the canonical base phi expansion of the natural numbers only vertical runs of 1's  with length a Lucas number occur, and all Lucas numbers occur as a run length. More precisely: in digit position $i$ only runs of length $L_{i-1}$ occur when $i\ge 1$, and only runs of length $L_{-i}$ occur when $i\le 0$.
\end{theorem}

\noindent  {\it Proof:}  The proof will be divided into five parts: $i=0$, $i=-1$, $i>0$ and $i<-1$ and $i$ even, $i<-1$ and $i$ odd.
 In the last three cases we partition the natural numbers in canonical Lucas intervals, and use the Recursive Structure Theorem  \ref{th:recCan}.
 For $i<-1$ the situation is more complicated than for $i>1$, which forces us to consider $i$ even and $i$ odd separately.

\medskip

\noindent {\bf Part\! 1:\! $\mathbf{i=0}$.}\; Then $L_{-i}=L_0=2$. We remark here that we did not mention in the statement of the theorem that for $i=0$ the first run deviates from the pattern: it has length 1, and this does not change if we would add $N=0$ to the table.

According to Theorem 5.1.~in \cite{Dekk-phi-FQ} one has $d_0(N)=1$   if and only if $N=\lfloor n\varphi\rfloor +2n +1$ for some natural number $n$. According to Lemma \ref{lem:gammabeta}, $\gamma(N)\ne \beta(N)$, if and only if there exists a natural number $n$, such that $N=\lfloor n\varphi \rfloor+2n$. If we combine these two statements we see that all the runs of 1's (except the first one) have length $2=L_0$ in digit position $i=0$.

\medskip

\noindent {\bf Part\! 2:\! $\mathbf{i=-1}$.\;}
 From Remark 5.1 from the paper \cite{Dekk-phi-FQ} we have that for the Bergman representation $d_{-1}(N)=1$  if  and only if $N=3\lfloor n\varphi\rfloor + n + 1$ for some natural number $n$. Since by Proposition \ref{prop: gammabeta} $\gamma(N)\ne \beta(N)$ if and only if $N$ is of type $\B$, which has $d_1d_0\cdot d_{-1}(N)=00\cdot 0$, we can deduce that also for the canonical representation $c_{-1}(N)=1$  if  and only if $N=3\lfloor n\varphi\rfloor + n + 1$ for some natural number $n$. This obviously implies that the runs of $1$'s at digit position $i=-1$ have length $L_1=1$.

\medskip

The proofs of Part 3, 4 and 5 are based on the Recursive Structure Theorem. When we perform the induction, we have to prove that runs of $1$'s do not extend beyond the intervals that are produced by the induction. In the following we will show that this holds for all digit positions with exception of $i=0$ and certain positions at the left end and right end of the expansion.  Note that we can ignore the case $i=0$, as it has already been dealt with in Part 1 of the proof.
For canonical Lucas intervals $\Gamma_{2n+2}$ this ``isolated run property" is considered in [*], for canonical Lucas intervals $\Gamma_{2n+1}$ in [**].

\medskip

[*] For the Recursive Structure Theorem, Part I,  we use that the runs of $1$'s in column $i$ in the interval $\Gamma_{2n+2}$ are a copy of the runs of $1$'s in the interval $[1,L_{2n+1}]$, \emph{except} for the left most column, which corresponds to digit position $L=2n+2$, and the right most column, corresponding to digit position $R=-2n-2$.

We still have to check that no new runs are created at the beginning or the end of the interval $\Gamma_{2n+2}$. This is obvious for the beginning, since $c_{i}(1)=0$ for $i\ne 0$, and for the end it follows from Lemma \ref{lem:ortho} [EO].

\medskip

[**]  For the Recursive Structure Theorem, Part II, we use that the runs of $1$'s in column $i$ in the interval $\Gamma_{2n+1}$ are a copy of the column of digit $i$ lying in the intervals $I_n, J_n$ and $K_n$, except for the three leftmost columns, and the five right most columns. These exceptions correspond to digit positions $L=2n+1, L-1=2n$, and $L-3=2n-1$, at the left, and positions except for the columns with indices $R+3=-2n+1, R+2=-2n, R+1=-2n-1, R=-2n-2$ when $N$ is from $I_n$ or $K_n$, and  except for the columns with indices $R+4, R+3, R+2, R+1, R$ when $N$ is from $J_n$. Here you use computations like $1000(10)^{-1}=100 1^{-1}$.

  This time,  we still have to check that no new runs are created at the beginning or the end of the intervals $I_n, J_n$ and $K_n$.

  For the beginning and the end of $\Gamma_{2n+1}=I_n\cup J_n\cup K_n$ this follows again from Lemma \ref{lem:ortho} [EO], respectively Lemma \ref{lem:ortho} [OE], except  for  $i= -2n-2$.

  By Remark \ref{rem:shifts} the digits in column $i$ of $I_n, J_n$ and $K_n$ are equal to the corresponding digits in column $i$ of the intervals $\Gamma_{2n-1},\Gamma_{2n-2}$ and $\Gamma_{2n-1}$.

 There will be no new run created on the boundary between $I_n$ and $J_n$. We have that $\Gamma_{2n-1}$ (the shift of $I_n$) ends with $N=L_{2n}$, and $c_i(L_{2n})=0$ by Equation (\ref{eq:gL1}), when  $c_i$ is not the digit of one of the last four columns.

 Finally, there will be no new run created on the boundary between $J_n$ and $K_n$, which are  shifts of the two successive intervals $\Gamma_{2n-2}$ and $\Gamma_{2n-1}$, by Lemma \ref{lem:ortho} [EO].

\bigskip

\noindent {\bf Part\! 3:\! {$\mathbf{i\ge 1}$.}}

We first illustrate how this works for the case $i=1$. The column of digit position 1 starts with a run of length $L_{i-1}=L_0=2$ in $\Gamma_1=[2,3].$
Then a 0 follows in $\Gamma_2=\{4\}$, followed by another run of length 2 in $\Gamma_3=[5,6,7]$. Suppose one has proved that only runs of length 2 occur in the Lucas intervals $\Gamma_1,\dots, \Gamma_m$ for some natural number $m$. We then proceed by induction, distinguishing the cases $m=2n$ and $m=2n+1$.

We start with the case $m=2n+1$. Then the next interval is $\Gamma_{2n+2}$. By the Recursive Structure Theorem part I, the column of digit position 1 lying in this interval is a copy of the column of digit 1 lying in the interval $[1,L_{2n+1}]$. Therefore, by the induction hypothesis, and [**] there will be only runs of length 2 in this part of the column.

For the case $m=2n$, the next interval is $\Gamma_{2n+1}$. By the Recursive Structure Theorem part II, the column of digit position 1 lying in this interval is a copy of the column of digit 1 lying in the intervals $I_n, J_n$ and $K_n$. Therefore, by the induction hypothesis and [*], there will be only runs of length 2 in this part of the column.  This ends the proof of the case $i=1$.

\medskip

We next consider the case $i$ for arbitrary $i\ge 2$. The proof is similar to the proof of the case $i=1$. The main complication is the change in the digits occurring at the left most part of the expansion  in the Recursive Structure Theorem.
This is solved by giving the induction more attention at the start.

The first run of $1$'s in digit column $i$ starts at the number $N=L_i+1$ in $\Gamma_i$, since all $c_i(N)=1$ for $N \in \Gamma_i$ (see $i=2n$ and $i=2n+1$ in Proposition \ref{prop:length}). This run has length  $L_{i-1}$, since $|\Gamma_i|=L_{i-1}$, and by Lemma \ref{lem:ortho}.

Next, for all $N$ in $\Gamma_{i+1}$ one has $c_i(N)=0$, simply because  $c_{i+1}(N)=1$. So no runs of $1$'s  occur in $\Gamma_{i+1}$.

We then pass to $\Gamma_{i+2}$.

Suppose $i$ is \underline{even}. Then, by the Recursive Structure Theorem part I, the column of digit position $i$ lying in the interval
$\Gamma_{i+2}$ is a copy of the column of digit $i$ lying in the interval $[1,L_{i+1}]$.
But we know already that this part only contains runs of $1$'s of length $L_{i-1}$ (actually there is a single run, lying in $\Gamma_i$.)
Then also $\Gamma_{i+2}$ will only have a run of $1$'s of length $L_{i-1}$. Moreover, this \emph{is} the length of that run by Lemma \ref{lem:ortho} [EO].

Suppose $i$ is \underline{odd}. Then, by the Recursive Structure Theorem part II  the column of digit position $i$ lying in the interval $\Gamma_{i+2}$ can be obtained from the column of digit $i$ lying in the intervals $\Gamma_{i},\Gamma_{i-1}$, and $\Gamma_{i}$. In the first case all the $L_{i-1}$ $1$'s turn into $0$'s, in the second case we obtain only $0$'s, simply because the right most $1$ of the expansions in $\Gamma_{i-1}$ is in column $i-1$, and in the third case  all the $L_{i-1}$  $1$'s turn into $1$'s.  Moreover, in this last case this \emph{is} a run of length $L_{i-1}$, by an application of Lemma \ref{lem:ortho} [OE]. Conclusion:  also $\Gamma_{i+2}$ will only have a run of $1$'s of length $L_{i-1}$.

\medskip

 Suppose one has proved that only runs of length $L_{i-1}$ occur in the Lucas intervals $\Gamma_i, \Gamma_{i+1},\dots, \Gamma_m$ for some natural number $m\ge i+2$. We then proceed by induction, distinguishing again the cases $m=2n$ and $m=2n+1$.

We start with the case $m=2n+1$. Then the next interval is $\Gamma_{2n+2}$. By the Recursive Structure Theorem part I, the column of digit position $i$ lying in this interval is a copy of the column of digit $i$ lying in the interval $[1,L_{2n+1}]=[1,L_m]$, except for the leftmost column. This column has digit position $2n+2$. But $i+2\le m=2n+1$,  therefore, by the induction hypothesis and [*], there will be only runs of length $L_{i-1}$ in this part of column $i$.

For the case $m=2n$, the next interval is $\Gamma_{2n+1}$. By the Recursive Structure Theorem part II, the column of digit position $i>0$ lying in this interval is a copy of the column of digit $i$ lying in the intervals $I_n, J_n$ and $K_n$, except for the three leftmost columns.
These three have indices $L=2n+1, L-1=2n$, and $L-3=2n-1$.
But $2n=m\ge i+2$, i.e., $i\le 2n-2$,  therefore, by the induction hypothesis and [**],  there will be only runs of length $L_{i-1}$ in this part of the column.

\medskip

\noindent {\bf Part\! 4:\! $\mathbf{i<-1}, \mathbf{i}$ even.\;}

Suppose $-i=2j$ is even. From Proposition \ref{prop:length} we obtain that the first run of numbers $N$ with digit $c_{-2j}$ equal to 1 starts with all numbers $N$ in the interval $\Gamma_{2j-1}$, and then \emph{continues} in the interval  $\Gamma_{2j}$. The run will not continue in the next interval $\Gamma_{2j+1}$, by Lemma \ref{lem:ortho}. Conclusion: the first run of $1$'s in column $-2j$ has length
$$|\Gamma_{2j-1}|+|\Gamma_{2j}|=L_{2j-2}+L_{2j-1}=L_{2j}=L_{-i}.$$

Next, we consider the interval $\Gamma_{2j+1}$. By the Recursive Structure Theorem part II, since both for $I_n$ and $K_n$ the last three digits
in the replacement equation are 001, there will be $0$'s in the corresponding parts of column $i$. The same is true for the part corresponding to the interval $J_n$.  So no runs of 1's occur in  $\Gamma_{2j+1}$. We then pass to $\Gamma_{2j+2}$.

Here, by the Recursive Structure Theorem part I, the column of digit position $i$ lying in the interval
$\Gamma_{2j+2}$ is a copy of the column of digit $i$ lying in the interval $[1,L_{2j+1}]$.
But we know already that this part only contains runs of $1$'s of length $L_{-i}$
(actually there is a single run, lying in the union of $\Gamma_{2j-1}$ and $\Gamma_{2j}$.)
Then also $\Gamma_{2j+2}$ will only have a run of $1$'s of length $L_{-i}$. Moreover, this \emph{is} the length of that run by Lemma \ref{lem:ortho} [EO].

Suppose one has proved that only runs of length $L_{-i}$ occur in the Lucas intervals $\Gamma_{-i+2}, \dots, \Gamma_{m}$ for some natural number  $m\ge 2j+2$. We then proceed by induction, distinguishing again the cases $m=2n$ and $m=2n+1$.

We start with the case $m=2n+1$. Then the next interval is $\Gamma_{2n+2}$. By the Recursive Structure Theorem part I, the column of digit position $i$ lying in this interval is a copy of the column of digit $i$ lying in the interval $[1,L_{2n+1}]$, except for the rightmost column. This column has digit position $-i=2n+2$. But $-i=2j\le m-2=2n-1$,  therefore, by the induction hypothesis and [**], there will be only runs of length $L_{-i}$ in this part of the column.

 For the case $m=2n$, the next interval is $\Gamma_{2n+1}$. By the Recursive Structure Theorem part II, the column of digit position $i$ lying in this interval is a copy of the column of digit $i$ lying in the intervals $I_n, J_n$ and $K_n$, except for the four rightmost columns for $I_n$ and $K_n$, and the five rightmost columns for $J_n$. But $-i\le R-5=2n-3$, since  $2n=m\ge 2j+2=-i+2$. Therefore, by the induction hypothesis, and [*] there will be only runs of length $L_{-i}$ in this part of the column. Here we need   $i\ne -(2n+2)$. This is satisfied because in fact, $-i<2n$.

\medskip

\noindent {\bf Part\! 5:\! $\mathbf{i<-1, \: i\: odd}$.}\;

Suppose  $-i=2j+1$ is odd. Where does the first run of $1$'s at digit position $i=-2j-1$ occur? This is more complicated than in all previous cases where this happened at position $L$ (for $i>0$) or position $R$ (for $i<0, i$ even).

\noindent Claim: the first run of $1$'s occurs in column $R+3$ in $\Gamma_{2j+3}$, where all numbers have $R=-2j-4$, as given in Proposition \ref{prop:length}.

Indeed, note that $\Gamma_{2j+1}$ and $\Gamma_{2j+2}$ would be the first two candidates for the occurrence of $1$'s at position $-(2j+1)$, but that both intervals have numbers $N$ with $R(N)=-(2j+2)$, and so there will be $0$'s at position $-(2j+1)$, since $11$ does not occur. The next candidate is the interval $\Gamma_{2j+3}$.
Here we use Lemma \ref{lem:R1001}, Part b), with $n=j+1$. This lemma gives that
$c_{-2j-1}(N)=1$ for the first $L_{2j+1}=L_{-i}$ numbers $N$ from $\Gamma_{2j+3}$, and $c_{-2j-1}(N)=0$ for the remaining numbers $N$.
This proves the claim above.

Next, consider the interval $\Gamma_{2j+4}$. By the Recursive Structure Theorem part I, the column of digit position $i$ lying in this interval is a copy of the column of digit $i$ lying in the interval $[1,L_{2j+3}]=\Gamma_0\cup\dots\cup\Gamma_{2j+2}$ (except for the rightmost column). But the first run of $1$'s occurs in column $R+3$ in $\Gamma_{2j+3}$, so there are no  $1$'s at all in column $R+3$ in $\Gamma_{2j+4}$.

\medskip

Next, suppose one has proved that only runs of length $L_{-i}$ occur in the canonical Lucas intervals $\Gamma_{-i+2}, \Gamma_{-i+1},\dots, \Gamma_{m}$ for some natural number $m\ge 2j+4$.  We then proceed by induction, distinguishing again the cases $m=2n$ and $m=2n+1$.

We start with the case $m=2n+1$. Then the next interval is $\Gamma_{2n+2}$. By the Recursive Structure Theorem part I, the column of digit position $i$ lying in this interval is a copy of the column of digit $i$ lying in the interval $[1,L_{2n+1}]=[1,L_m]$, except for the rightmost column with index $R=-2n-2$. So we need that $-i<-R=2n+2$, which holds iff $2j+1<m$. This is certainly satisfied.

Therefore, by the induction hypothesis and [*], there will be only runs of length $L_{-i}$ in this part of the column of digit $i$.

 For the case $m=2n$, the next interval is $\Gamma_{2n+1}$, with $R=-2n-2$. By the Recursive Structure Theorem part II, the column of digit position $i$ lying in this interval is a copy of the column of digit position $i$ lying in the intervals $I_n, J_n$ and $K_n$, except for the columns with indices $R+3, R+2, R+1, R$ when $N$ is from $I_n$ or $K_n$, and  except for the columns with indices $R+4, R+3, R+2, R+1, R$ when $N$ is from $J_n$.
So we need that $-i<-R-4=2n+2-4=m-2$, which holds iff $2j+1<m-2$, which is satisfied since $m\ge 2j+4$.

 Therefore, by the induction hypothesis and [**], there will be only runs of length $L_{-i}$ in this part of the column.  \hfill $\Box$

\section{Final remarks}

\subsection{Two-dimensional characteristic}\label{sec:third}

There is a third characteristic of representations of the natural numbers, which might be labelled as `two-dimensional'. This amounts to the observation that the lengths of the runs of $1$'s spread over the table of expansions in chains of consecutive Lucas numbers, cf.~Figure \ref{Fig:chains}.
These chains are finite, except the first one, which starts at $N=2$, and consists of the runs of the  $c_L$ digits.

The pattern consists of two kinds of chains:
\begin{enumerate}
	\item	At the left side, i.e., $i>0$,  the lengths of the links in the chains follow $(L_{n})_{n\geq0}= 2,1,3,4,7,11,\dots$.
	\item	At the right side, i.e., $i\le 0$,  the lengths of the links in the chains follow $(L_{n})_{n\leq 0}= 2,-1,3,-4,7,-11,\dots$.
\end{enumerate}
Here the sign of the length indicates the direction in which the link goes.

\begin{figure}[h]%
\begin{center}
{\footnotesize
\begin{tabular}{|p{6mm}|p{3mm}|p{3mm}|p{2mm}|p{2mm}|p{2mm}|@{}p{0.3mm}@{}|p{2mm}|p{4mm}|p{4mm}|p{4mm}|p{4mm}|p{4mm}|p{6mm}|}
\!\!{$\mathbf (L_{n})$}  & $c_5$  & $c_4$&$c_3$  &$c_2$&$c_1$ & \!\! &$c_0$ &$c_{-1}$&$c_{-2}$&$c_{-3}$&$c_{-4}$&$c_{-5}$&$(L_{n})$\\
	\hline
	2$\blacktriangledown$ &  &  &  &  & 1 &   &   &  &  &  &  1 &  &   7$\blacktriangledown$ \\
	   &  &  &  &  & 1 &   &   &  &  &   & 1 &  &  \\
	1$\blacktriangledown$ &  &  &  & 1 &  &   &   &  &  &   & 1 &  &   \\
	3$\blacktriangledown$  &  &  & 1 &  &  &   &   &  &  &   & 1 &  &   \\
	   &  &  & 1 &  &  &   &   &  & 1 &  & 1 &  &     3$\blacktriangledown$ \\
	   &  &  & 1 &  &  &   & 1 &  & 1 &  & 1 &  &      2$\blacktriangledown$ \\
	4$\blacktriangledown$  &  & 1 &  &  &  &   & 1 & &  1 &  &  1 & &  \\
	   &  & 1 &  &  &  &   &   & 1&  & 1 & & 1 &   -1$\blacktriangle$  \\
	   &  & 1 &  &  &  &   &   &  &  & 1 &  & 1 &  \\
	   &  & 1 &  &  &  &   &   &  &  & 1 &  & 1 &  \\
	7$\blacktriangledown$  & 1 &  &  &  &  &   &   &  &  & 1 &  & 1 &    -4$\blacktriangle$\\
	     & 1 &  &  &  &  &   &   &  &  &  &  & 1 &  \\
	     & 1 &  &  &  &  &   &   &  &  &  &  & 1 &  \\
	     & 1 &  &  &  &  &   &   &  &  &  &  & 1 &  \\
	     & 1 &  &  &  &  &   &   &  &  &  &  & 1 &  \\
	     & 1 &  &  &  &  &   &   &  &  &  &  & 1 &  \\
	     & 1 &  &  &  &  &   &   &  &  &  &  & 1 &  \\	
	     &   &  &  &  &  &   &   &  &  &  &  & 1 &   \!-11$\blacktriangle$\!
         \end{tabular}}
 \end{center}
 \caption{\footnotesize Two typical chains. Leftmost column: $L_n, n\ge 0$, rightmost column: $L_n, n\le 0$.}\label{Fig:chains}
\end{figure}


\subsection{Positions}
We conjecture that the positions at which the $1$'s occur in the column of digit $i$  are given by unions of generalized Beatty sequences. Generalized Beatty sequences, defined in \cite{GBS}, are sequences $V(p,q,r)=(V_n)$ of the form $V_n = p\lfloor n \alpha \rfloor + q\,n +r $, $n\ge 1$,  where $\alpha$ is a real number, and $p,q,$ and $r$ are integers. Note that this has been proved in our paper for $i=0$ and $i=-1$.

\subsection{Generalizations}  We believe that for other irrational numbers than the golden mean our approach makes sense. In particular for the metallic means, a special class of numbers contained in the class given in Theorem 2 of \cite{FrouSol}. We shortly discuss the case of the silver mean $\sigma:=1+\sqrt{2}$.

The standard representation of the natural numbers in base $\sigma$ is given by 
 $$N= \sum_{i=-\infty}^{\infty} d_i \sigma^i,\vspace*{-.0cm}$$
  with digits $d_i=0,1$ or 2, and where $d_{i+1}d_i = 21$ or $22$ is not allowed.
  
  The role of the Lucas numbers $(L_n)$ is now taken over by the Pell-Lucas numbers $(\LP_n)=2,2,6,14,34,\dots$, defined by
  $$\LP_0=2, \quad\LP_1=2,\quad  \LP_{n+2}= 2\LP_{n+1}+ \LP_n,\quad{\rm for}\;n=0,1,2,\dots $$

 We write $\BL(N)$ for the standard expansion of $N$ in base $\sigma$, and  $\GL(N)$ for the canonical expansion of $N$ in base $\sigma$.
 This time canonical means  that the digits are  $c_i=0,1$ or 2, and that $c_{i+1}c_i = 21$ or $22$ is not allowed, \emph{except} that $c_1c_0 = 21$, as soon as this is possible.

 The following table displays these representations.

 \bigskip

 \begin{tabular}{|r|c|c|}
\hline
  \; $N^{\phantom{|}}$ & $\BL(N)$ &\; $\GL(N)$ \\[.0cm]
\hline
1  \; & \;\,1$\cdot$0         &  \:1$\cdot$0               \\
2  \; & \;\,2$\cdot$0         &  \:2$\cdot$0            \\
3  \; &   10$\cdot$11          &  10$\cdot$11            \\
4  \; &   11$\cdot$11       &  11$\cdot$11           \\
5  \; &   20$\cdot$01     &  20$\cdot$01         \\
6  \; & \! 100$\cdot$01    &  21$\cdot$01         \\
\hline
\end{tabular} \qquad
\begin{tabular}{|r|c|c|}
\hline
  \; $N^{\phantom{|}}$ & $\BL(N)$ &\; $\GL(N)$ \\[.0cm]
\hline
7  \; & 101$\cdot$01    &    101$\cdot$01        \\
8  \; & 102$\cdot$01    &    102$\cdot$01        \\
9  \; & 110$\cdot$12   &     110$\cdot$12       \\
10 \; & 111$\cdot$12    &    111$\cdot$12        \\
11 \; & 120$\cdot$02     &   120$\cdot$02         \\
12 \; & 200$\cdot$02  &      121$\cdot$02     \\
\hline
\end{tabular} \qquad
\begin{tabular}{|r|c|c|}
\hline
  \; $N^{\phantom{|}}$ & $\BL(N)$ &\; $\GL(N)$ \\[.0cm]
\hline
13  \; & 201$\cdot$02     &    201$\cdot$02        \\
14 \; & 202$\cdot$02      &   202$\cdot$02         \\
15 \; &\: 1000$\cdot$2011   &  \: 1000$\cdot$2011      \\
16 \; &\: 1001$\cdot$2011   &  \: 1001$\cdot$2011      \\
17 \; &\: 1010$\cdot$1011   & \:  1010$\cdot$1011     \\
18 \; &\: 1011$\cdot$1011   & \:  1011$\cdot$1011     \\
\hline
\end{tabular}

\bigskip

We conjecture that the sequence of natural numbers $6,12,20,26,34,\dots$ for which  $\BL(N)\ne \GL(N)$ is equal to the generalized Beatty sequence (with $\alpha=\sigma$)
$V(2,2,0)=(2\lfloor n(\sigma+1) \rfloor)$.

\medskip

What are the intervals of constant expansion length for the two representations by powers of $\sigma$?\\
In the same way as in Section \ref{sec:length} we define the   \emph{Pell-Lucas intervals} and the \emph{canonical Pell-Lucas intervals}  :
\begin{eqnarray*}
  \LambdaL_0:=\{1,2\},&\,& \LambdaL_{2n}:=[\LP_{2n},\LP_{2n+1}]\;{\rm for}\:n\ge 1,\quad \LambdaL_{2n+1}:=[\LP_{2n+1}+1, \LP_{2n+2}-1]\;{\rm for}\:n\ge 0 \\
 \GammaL_0:=\{1,2\},&\,& \GammaL_{n}:=[\LP_{n}+1,\LP_{n+1}]\;\;{\rm for}\:n\ge 1
\end{eqnarray*}
We conjecture that the $\LambdaL_n$ are the intervals of constant expansion length of the standard silver mean representation, and that the $\GammaL_n$ are the intervals of constant  expansion length  of the canonical silver mean representation.

\medskip

We next consider vertical runs. Here there are runs of $1$'s and runs $2$'s. We conjecture that in the column of digit $i>0$ there are only runs of $1$'s of length $\LP_{i}$,
followed directly by runs of $2$'s of length $\LP_{i-1}$.

We also conjecture that for $i < 0$ there are either runs of 1's of length $\LP_{-i}$ or runs of 1's of length $\LP_{-i}$ directly followed by runs of 2's of length $\LP_{-i}$. In the odd columns the order of the runs of 1's and runs of 2's is reversed. So, for example, for $i=-1$  between $0$'s, there are only blocks 11 and 2211. Note that these reversals are in line with the changes of direction of the base phi expansions observed at the end of Section \ref{sec:third}.

\end{document}